\documentclass[10pt]{article}
\usepackage[top=4.2cm, bottom=4.2cm, left=4.0cm, right=4.0cm]{geometry}
\usepackage{amssymb}
\usepackage{latexsym}
\usepackage{parskip}
\usepackage{amsthm}
\usepackage[all]{xy}
\usepackage[T1]{fontenc}

\newtheorem{thm}{Theorem} \newtheorem{prop}{Proposition}
\newtheorem{lemma}{Lemma} \newtheorem{cor}[thm]{Corollary}

\theoremstyle{definition}

 \title{{ An Example of non-quenched Convergence in the Conditional CLT for Discrete Fourier Transforms \\}} \author{{\Large David Barrera}\\{} \\ {Department of Mathematical Sciences} \\{  University of Cincinnati}\\ {\small PO Box 210025,
Cincinnati, Oh 45221-0025, USA.}\\{\small {\it Email}: barrerjd@mail.uc.edu}}

 \date{} 
 \begin{document} \maketitle

\begin{abstract}
A recent result by Barrera and Peligrad (\cite{BaPe}, Theorem 1) shows that the quenched Central Limit Theorem holds for the components of the discrete Fourier transforms of a stationary process in $L^2$ orthogonal to the subspace of functions that are measurable with respect to the initial sigma-field. In this note we address the question of whether the quenched CLT remains true for the Fourier transforms without taking orthogonal projections, as could be expected in view of previous,  related results about the annealed convergence of the process under consideration (see for instance \cite{PeWu}, Theorem 2.1). 

We give a negative answer to this question by exhibiting an example of a process satisfying the hypothesis of Theorem 1 in \cite{BaPe} for which the Fourier transforms do not satisfy a quenched limit theorem. The proof combines ideas from a construction due to Voln\'{y} and Woodroofe (see \cite{VoWo})  with an interpretation of the results in \cite{BaPe} in the case of linear processes, and with applications of some previous results related to discrete Fourier transforms.

\end{abstract}



\section{Introduction and Notation.}

Let $(X_{n})_{n\in\mathbb{Z}}$ be a strictly stationary centered, ergodic sequence of random variables defined on a probability space $(\Omega, \mathcal{F}, \mathbb{P})$. This is: $X_{n}=X_{0}\circ T^{n}$,  where $T:\Omega\to \Omega$ is an ergodic, bimeasurable, invertible transformation and $EX_{0}=0$. Assume that $X_{0}\in L^{2}_{\mathbb{P}}(\mathcal{F}_{0})$ where $\mathcal{F}_{0}\subset\mathcal{F}$ is a sigma algebra satisfying $\mathcal{F}_{0}\subset T^{-1}\mathcal{F}_{0}$ (i.e. $T^{-1}$ is $\mathcal{F}_{0}$ measurable), define $\mathcal{F}_{n}:=T^{-n}\mathcal{F}_{0}$ for all $n\in\mathbb{Z}$ and $\mathcal{F}_{-\infty}:=\cap_{n\in \mathbb{Z}}\mathcal{F}_{n}$. Denote by $E_{n}$ the conditional expectation with respect to $\mathcal{F}_{n}$. So $E_{n}Z:=E[Z|\mathcal{F}_{n}]$ for every integrable random variable $Z$.

Define, for every  $\theta\in[0,2\pi)$ the {\it $n-$th discrete Fourier transform of  $(X_{k}(\omega))_{k}$ at $\theta$} by
\begin{equation}
\label{nfoutradef}
S_{n}(\theta,\omega):=\sum_{k=0}^{n-1}e^{ik\theta}X_{k}({\omega}).
\end{equation}

When $\theta\in(0,2\pi)$ is fixed, we will denote by $S_{n}(\theta)$ the random variable $S_{n}(\theta,\cdot)$. In the special case $\theta=0$ we denote by $S_{n}$ the random variable $S_{n}(0,\cdot)$. So $S_{n}(\omega):=\sum_{k=0}^{n-1}X_{n}(\omega)$.

Assume also that $E_{0}$ is regular. This is, that there exists a family of measures $\{\mathbb{P}_{\omega}\}_{\omega\in \Omega}$ such that for every integrable function $X$, 
$$\omega\mapsto \int_{\Omega}X(\omega')d\mathbb{P}_{\omega}(\omega')$$
defines a version of $E_{0}X$.

Finally, denote by $\lambda$ the Lebesgue measure on $[0,2\pi)$ (or any other Borelian in $\mathbb{R}$).

Under these assumptions Barrera and Peligrad, in \cite{BaPe}, proved the following theorem:

\bigskip

\begin{thm}
\label{BaPeTh}
There exists a set $I\subset (0,2\pi)$ with $\lambda(I)=2\pi$ such that, for all $\theta\in I$, the complex-valued random variables

\begin{equation}
\label{norparsumc}
Y_{n}(\theta):=\frac{1}{\sqrt{n}}(S_{n}(\theta)-E_{0}S_{n}(\theta))
\end{equation}

converge to a complex Gaussian random variable under $\mathbb{P}_{\omega}$ for all $\omega$ in a set $\Omega_{\theta}$ with $\mathbb{P}(\Omega_{\theta})=1$. The asymptotic distribution of the real and imaginary parts corresponds to a bivariate Gaussian random variable with independent entries, each with mean zero and variance

$$\sigma_{\theta}^2=\lim_{n\to\infty}\frac{E_{0}|S_{n}(\theta)-E_{0}S_{n}(\theta)|^2}{2n}.$$

(the limit exists with probability one,  and it is nonrandom). 
\end{thm}

\section{Quenched Convergence}

In the context of the previous section, and given a distribution function $F_{Y}$ (associated to some random variable $Y$ defined on a probability space $(\Omega',\mathcal{F}',\mathbb{P}')$), we say that  the process {\it $Y_{n}$ converges to $Y$ in the quenched sense}, denoted here by $Y_{n}\Rightarrow_{Q} Y$, if for almost every $\omega$, and every continuous and bounded function $f$, $E^{\omega}f(Y_{n})\to_{n\to \infty} Ef(Y)$, 
where $E^{\omega}$ denotes integration with respect to $\mathbb{P}_{\omega}$ and $Ef(Y):=\int_{\mathbb{R}}f(y)\, dF_{Y}(y)=\int_{\Omega'}f\circ Y(\omega')\,d\mathbb{P}'(\omega')$.

In other words, we require that $E_{0}f(Y_{n})\to Ef(Y)$ a.s. (over a set not depending  of $f$). By Portmanteau's theorem, this amounts to say that for almost every $\omega$, $\mathbb{P}_{\omega}(Y_{n}\leq y)\to F_{Y}(y)$ at every continuity point $y$ of $F_{Y}$. Theorem \ref{BaPeTh} is thus a statement about the quenched convergence of $Y_{n}(\theta)$ for every $\theta\in I$.

Note that quenched convergence implies convergence in distribution (``annealed'' convergence) by the dominated convergence theorem: for every uniformly continuous and bounded function $f$,
$$\int_{\Omega}f(Y_{n})\, d\mathbb{P}(\omega)=\int_{\Omega}E_{0}f(Y_{n})(\omega)\, d\mathbb{P}(\omega)\to_{n\to \infty} \int_{\Omega} Ef(Y)\, d\mathbb{P}(\omega)=Ef(Y).$$

In particular, Theorem \ref{BaPeTh} relates to some previous results about annealed convergence (see for instance \cite{PeWu} and the references therein).

It is worth to remark that, without additional assumptions, the methods of \cite{BaPe} do not give a description of the elements in the set $I$. The martingale version of the theorem, used to approximate the general case, works provided that $e^{-2it}$ is not an eigenvalue of the Koopman operator associated to $T$ (namely $f\mapsto f\circ T$ for all $f\in L^2_{\mathbb{P}}(\Omega,\mathbb{C})$), and therefore we consider these as exceptional values. A consideration of the classical case $\theta=0$ shows that more hypotheses may be needed to guarantee the convergence in distribution of $Y_{n}(\theta)$ outside of $I$.

\subsection{Possible limit Laws for a given initial point.}

Suppose that we know of an integrable process $(Y_{n})_{n}$ that $Y_{n}-E_{0}Y_{n}\Rightarrow_{Q} Y$, say $Y_{n}-E_{0}Y_{n}\Rightarrow Y$ under $\mathbb{P}_{\omega}$ for $\omega\in \Omega_{0}$, where $\mathbb{P}(\Omega_{0})=1$. What are the possible limit laws for $Y$ under $\mathbb{P}_{\omega}$, for a fixed $\omega$?

To answer this question we depart from the following result (see the proof of Lemma 18 in \cite{BaPe}): 

\medskip

\begin{lemma}
\label{BarPelLem18}

If $\mathcal{F}_{0}\subset\mathcal{F}$ is a $\sigma-$algebra for which $E[\,\cdot\,|\mathcal{F}_{0}]=:E_{0}$ admits a regular version in the sense explained above {(see the introduction)}, $Y$ is $\mathcal{F}_{0}-$measurable, and $X$ is a given random variable, then there exists $\Omega_{1,Y}\subset \Omega$ with $\mathbb{P}(\Omega_{Y})=1$ such that, for every $g:\mathbb{R}^{2}\to \mathbb{R}$  continuous and bounded
\begin{equation}
\label{lecoex}
E^{\omega}[g(X,Y(\omega))]=E^{\omega}[g(X,Y)]
\end{equation}
for all $\omega\in \Omega_{Y}$.
\end{lemma}

This lemma, in combination with Proposition \ref{quecon} and Lemma \ref{degcascon} in the appendix, gives the following proposition.

\medskip

\begin{prop}[Possible limit Laws for a fixed starting point.]
\label{queconrev}
With the notation of Lemma \ref{BarPelLem18}, assume that $(Y_{n})_{n}$ is an integrable process such that $Y_{n}-E_{0}Y_{n} \Rightarrow Y$ under $\mathbb{P}_{\omega}$ for all $\omega\in \Omega_{0}$ ($\Omega_{0}\subset \Omega$ is {\it any} given set, not even assumed measurable), and let $\Omega_{1}:=\cap_{n}\Omega_{Y_{n}}$. Then given $\omega\in \Omega_{0}\cap\Omega_{1}$, $Y_{n}$ is convergent under $\mathbb{P}_{\omega}$ to some random variable $Z_{\omega}$ if and only if $L(\omega)=\lim_{n\to\infty}E_{0}Y_{n}(\omega)$ exists, in which case $Z_{\omega}=Y+ L(\omega)$ (in distribution).
\end{prop}

\medskip

%{\bf Remark:} A modification of the proof of the necessity part of Lemma \ref{BaPeTh2}, using Proposition \ref{quecon} below, reveals the following: in the context of Theorem \ref{BaPeTh}, and given $\theta \in I$ for which $\sigma_{\theta}^{2}>0$, if there exists a non constant limit $Z_{\omega}(\theta)$ for $Z_{n}(\theta)$ under $\mathbb{P}_{\omega}$ %for $\mathbb{P}$-a.e $\omega$ (resp: 
%for $\omega$ in a set of probability $\eta>0$, then $L(\omega)=\lim_{n}E_{0}S_{n}(\theta)(\omega)/\sqrt{n}$ is well defined
%$\mathbb{P}$-a.s (resp: 
%with probability $\eta$ (this is, the limit exists and it is finite
% $\mathbb{P}$-a.s, 
%over a set of probability $\eta$), and if we denote by $Y(\theta)$ the random complex function whose distribution describes the asymptotic limit of $Y_{n}(\theta)$, then for $\omega$ in a set of probability $\eta$, $Z_{\omega}(\theta)= Y(\theta)+L(\omega)$ (in distribution). In particular, it is not possible to have a nondegenerate asymptotic limit for $Z_{n}(\theta)$ under $\mathbb{P}_{\omega}$ for $\omega$ in a nonnegligible set if one is able to prove that, with probability one, $\limsup_{n}|E_{0}S_{n}(\theta)|/\sqrt{n}=\infty$. In the case $\sigma_{\theta}^{2}=0$ a similar argument (using part (b) of Proposition \ref{quecon}) gives that, if $\omega\in \Omega_{\theta}$ and $\limsup_{n}|E_{0}S_{n}(\theta)(\omega)|/\sqrt{n}=\infty$, then $Z_{n}(\theta)$ cannot converge under $\mathbb{P}_{\omega}$ (except perhaps for $\omega$ in a set of $\mathbb{P}-$measure zero).

{\bf Proof:} Given $\omega\in \Omega_{0}\cap\Omega_{1}$ and any bounded and continuous function $g$
$$E^{\omega}g(Y_{n}-E_{0}Y_{n}(\omega))=E^{\omega}g(Y_{n}-E_{0}Y_{n})\to_{n} Eg(Y),$$
so that  $Y_{n}-E_{0}Y_{n}(\omega)\Rightarrow Y$ under $\mathbb{P}_{\omega}$. From $Y_{n}=Y_{n}-E_{0}Y_{n}(\omega)+E_{0}Y_{n}(\omega)$ the conclusion follows via Proposition \ref{quecon} if $Y$ is not constant and via Lemma \ref{degcascon} if $Y$ is constant.\qed

\subsection{The Question} 

Let us define, for every $\theta\in[0,2\pi)$, $Z_{n}(\theta):=\frac{1}{\sqrt{n}}S_{n}(\theta)$. It is natural to ask whether the random centering in (\ref{norparsumc}) is necessary to obtain quenched convergence. For regular processes (namely $E[X_{0}|\mathcal{F}_{-\infty}]=0$) this amounts, in view of Theorem 2.1 in \cite{PeWu}, to whether the conclusion of Theorem 1 holds with $Z_{n}(\theta)$ in place of $Y_{n}(\theta)$. 

The first observation in this direction is given by the following lemma.

\bigskip

\begin{lemma}
\label{BaPeTh2}
With the notation of Theorem \ref{BaPeTh}, fix $\theta\in (0,2\pi)$ ($\theta$ may or may not be in $I$), and assume that $Y_{n}(\theta)\Rightarrow_{Q} Y_{\theta}$ for some $Y_{\theta}$. Then $Z_{n}(\theta)\Rightarrow_{Q} Y_{\theta}$ if and only if $E_{0}S_{n}(\theta)=o(\sqrt{n})$ almost surely. 
\end{lemma}

{\bf Proof:} {\itshape Sufficiency.} First note that for any set $A\in \mathcal{F}$, $\mathbb{P}(A)=\int_{A}\mathbb{P}_{\omega}(A)d\,\mathbb{P}({\omega})$. Applying this observation to the complement of $[E_{0}S_{n}(\theta)=o(\sqrt{n})]$ we see that if $E_{0}S_{n}(\theta)=o(\sqrt{n})$ a.s. then $E_{0}S_{n}(\theta)=o(\sqrt{n})$  $\mathbb{P}_{\omega}-$a.s. for $\mathbb{P}-$a.e. $\omega$, and therefore the asymptotic distributions of $Y_{n}(\theta)$ and $Z_{n}(\theta)$ (if any) must be the same under $\mathbb{P}_{\omega}$ for $\mathbb{P}-$ a.e $\omega$. This proves sufficiency.

{\itshape  Necessity.}  We  appeal to Proposition \ref{queconrev} (applied to the real and imaginary parts of the processes in question) by taking, in place of $\Omega_{0}$, $\Omega_{\theta}:=\{\omega\in \Omega: Y_{n}(\theta)\Rightarrow Y_{\theta} \mbox{\, \,\,\it under  $\mathbb{P}_{\omega}$}\}$. This gives that $Y_{\theta}=Y_{\theta}+\lim_{n}{E_{0}S_{n}(\theta,\omega)}/{\sqrt{n}}$ for $\omega\in \Omega_{1,\theta}:=\Omega_{1}\cap\Omega_{\theta}$. This is, that 
$$\lim_{n}\frac{E_{0}S_{n}(\theta,\omega)}{\sqrt{n}}=0  \mbox{\,\,\,\,\, \it for all $\omega\in \Omega_{1,\theta}$ }. $$\qed

Therefore, to give a negative answer to our question, we must provide a regular process $(X_{n})_{n}$ satisfying the hypothesis of Theorem \ref{BaPeTh} for which
\begin{equation}
\label{protoprov}
\mathbb{P}\left(\limsup_{n\to \infty}\left|\frac{1}{\sqrt{n}}E_{0}S_{n}(\theta)\right|>0\right)>0\mbox{\,\,\,\, \it for   $\theta$ in a set $I'$ with $\lambda(I')>0$.} 
\end{equation}

Proving (\ref{protoprov}) gives, in particular, the necessity of random centering for a nonnegligible subset of $I$ (namely $I\cap I'$).

In their paper \cite{VoWo}, Volny and Woodroofe provide an example of a sequence $(X_{n})_{n}$ for which a quenched CLT holds for $(Y_{n}(0))_{n}$ but not for $(Z_{n}(0))_{n}$. In this paper, we adapt their construction to give an example satisfying (\ref{protoprov}) with $I'=[0,2\pi)$. The random centering ``$-E_{0}S_{n}(\theta)$'' is therefore a necessary condition for the conclusion of Theorem \ref{BaPeTh} to hold. 

The main novelty adapting the example in \cite{VoWo}, which arises from a careful construction of a sequence $(a_{n})_{n}$ of nonnegative coefficients of a linear process is that, in order to guarantee the validity of the ``inductive step'' defining $a_{n+1}$ from $a_{1},\dots, a_{n}$, one needs to prove that a certain type of convergence is uniform in $\theta$ (see Lemma \ref{lemxij} below). While it would be sufficient to prove this uniform convergence for $\theta$ in an open subinterval $I'$ of $[0,2\pi)$ in order to construct a valid example, a compactness argument allows us to do it for $I'=[0,2\pi)$.

We will give an example of a process which, indeed, has the following property

\begin{equation}
\label{protoprovinf}
\mathbb{P}\left(\limsup_{n\to \infty}\left|\frac{1}{\sqrt{n}}E_{0}S_{n}(\theta)\right|=\infty\right)=1\mbox{\,\,\,\, \it for   all $\theta\in [0,2\pi)$.} 
\end{equation}

For this process, Proposition \ref{queconrev} shows that $\frac{1}{\sqrt{n}}S_{n}(\theta)$ {cannot} admit an asymptotic limit under $\mathbb{P}_{\omega}$ for $\mathbb{P}$-a.e $\omega$.

The rest of the paper is presented as follows: in Section 3 we specialize our study to the case in which $(X_{k})_{k\in\mathbb{Z}}$ is a linear process generated by convolution of a sequence of i.i.d random variables and a sequence in $l^{2}(\mathbb{N})$. For this family of processes, we give an interpretation of the results above in terms of convergence of Fourier series of perodic functions (Proposition \ref{asyvart}), and introduce a result necessary to construct the example (Lemma \ref{lemxij}). In Section 4 we present the construction itself. 

\section{The Fourier Transforms of a Linear Process}

Let $(a_{n})_{n\in\mathbb{N}}$ be a sequence in $l^{2}(\mathbb{N})$ (namely $\sum_{n}a_{n}^2<\infty$), and let $(\xi_{k})_{k\in\mathbb{Z}}$ be an i.i.d. sequence of centered, square integrable random variables defined on some probability space  $(\Omega',\mathcal{F}',\mathbb{P}')$. The {\it linear process} $(X_{k})_{k\in \mathbb{Z}}$ generated by $(a_{k})$ and $(\xi_{k})$ is defined by

\begin{equation}
\label{linprodef}
X_{k}:=\sum_{j\in\mathbb{N}}a_{j}\,\xi_{k-j}.
\end{equation}

The orthogonality of $(\xi_{k})_{k}$ shows that $X_{k}$ is well defined as en element of $L^2_{\mathbb{P}\,'}$ and that $EX_{k}^{2}=||\xi_{0}||_{_2}^2\sum_{j}{a_{j}^2}$. 

We can regard $\Omega=\mathbb{R}^{\mathbb{Z}}$ as a probability space whose sigma algebra is the product sigma algebra and whose probability measure is $\mathbb{P}=\mathbb{P}'\xi^{-1}$, where $\xi:\Omega'\to \Omega$ is given by $\xi(\omega'):=(\xi_{j}(\omega'))_{j\in \mathbb{Z}}$. It is well known that, with this structure (because $(\xi_{k})_{k}$ is i.i.d.), the left shift $T:\Omega\to \Omega$ characterized by $x_{k}\circ T=x_{k+1}$, where $x_{j}:S\to \mathbb{R}$ is the projection on the $j-$th coordinate, is weakly mixing (and therefore also ergodic).

Note that the coordinate functions $x_{j}$ are a ``copy'' of the sequence $(\xi_{j})$: they are independent and have the same distribution. In particular, $X_{k}$ can be regarded as the function $X_{k}:\Omega\to \mathbb{R}$ given by $X_{k}(\omega)=X_{k}((x_{j}(\omega))_{j}):=\sum_{j}a_{j}x_{k-j}(\omega)$. 

Clearly, $X_{k}=X_{0}\circ T^{k}$. In this way $(X_{k})_{k\in\mathbb{Z}}$ can be interpreted as a strictly stationary, centered, and ergodic sequence in $L^{2}_{\mathbb{P}}$. 

In this case, we choose $\mathcal{F}_{n}:=\sigma((x_{k})_{k\leq n})$ for all $n\in \mathbb{Z}$ and we define $\mathbb{P}_{\omega}$ as the measure corresponding to ``partial integration with respect to the future''. This is: given $\omega_{0}\in \Omega$, $\mathbb{P}_{\omega_{0}}=\mathbb{P}\pi_{\omega_{0}}^{-1}$ where $\pi_{\omega_{0}}:S\to S$ is given by 

$$x_{k}(\pi_{\omega_{0}}(\omega))= \left\{ 
  \begin{array}{l l}
     x_{k}(\omega_{0})& \quad \textrm{if $k\leq 0$ } \\
    x_{k}(\omega) & \quad \textrm{if $k>0$}\\
    
  \end{array} \right .$$

This brings us to the hypohteses at the beginning. By Kolmogorov's zero-one law, $(X_{n})_{n}$ is a regular process.

As $T$ is weakly mixing its only eigenvalue is $\lambda_{0}=1$ (see for instance \cite{Pe}, p.65).

\subsection*{Coefficients of the Fourier Transforms}

Under the given hypothesis ($(a_{n})_{n}\in l^{2}(\mathbb{N})$), Carleson's theorem (\cite{Car}) guarantees the convergence a.s.  of the series

\begin{equation}
\label{detpart}
f({\theta})=\sum_{j\geq 0} a_{j}e^{ij\theta}            
\end{equation}
and $f(\theta)$, thus defined, is a $2\pi-$periodic function, square integrable over $[0,2\pi)$, and satisfying $\hat{f}(n)=a_{n}$, where $\hat{f}$ denotes the Fourier transform

$$\hat{f}(x)=\frac{1}{2\pi}\int_{0}^{2\pi}e^{-ixy}f(y)\,d\lambda(y).$$

Denote by $f_{k}(\theta):=\sum_{j=0}^{k-1}a_{j}e^{ij\theta}$ ($f_{k}=0$ if $k\leq 0$). Then we have the following two expressions for $S_{n}(\theta)$:

\begin{equation}
\label{snexp1}
S_{n}(\theta)=\sum_{j=-\infty}^{n-1}(f_{-j+n}-f_{-j})(\theta)\xi_{j}e^{ij\theta},
\end{equation}

\begin{equation}
\label{snexp2}
S_{n}(\theta)=\sum_{k=0}^{\infty}a_{k}\sum_{j=0}^{n-1}e^{ij\theta}\xi_{j-k}.
\end{equation}

Now let's denote, for all $k\geq 0$,  

\begin{equation}
\label{defzeta}
\zeta_{-k}(\theta):=\sum_{j=0}^{k}e^{-ij\theta}\xi_{-j} \mbox{ ($\zeta_{-k}=0$ if $k<0 $)}. 
\end{equation}

Then from  (\ref{snexp1}) and (\ref{snexp2}) the following two equalities follow respectively: 

\begin{equation}
\label{exprojsn1}
E_{0}S_{n}(\theta)=\sum_{j\leq 0}\xi_{j}(f_{-j+n}-f_{-j})(\theta)e^{ij\theta},
\end{equation}

\begin{equation}
\label{exprojsn2}
E_{0}S_{n}(\theta)=\sum_{j\geq 0} a_{j}(\zeta_{-j}-\zeta_{-j+n})(\theta)e^{ij\theta}.
\end{equation}

In particular
$$E_{0}|S_{n}(\theta)-E_{0}S_{n}(\theta)|^2=E_{0}|\sum_{j=1}^{n-1}e^{ij\theta}\xi_{j}f_{-j+n}(\theta)|^2=$$
$$||\xi_{0}||_{_2}^2\sum_{j=1}^{n-1}|f_{n-j}(\theta)|^2
=||\xi_{0}||_{_2}^2\sum_{j=1}^{n-1}|f_{j}(\theta)|^2,$$
so that 
$$\lim_{n\to \infty}\frac{E_{0}|S_{n}(\theta)-E_{0}S_{n}(\theta)|^2}{n}=\lim_{n\to \infty}\frac{1}{n}\sum_{j=1}^{n-1}|\,||\xi_0||_2 f_{j}(\theta)|^2=|\,||\xi_0||_{_2} f(\theta)|^2.$$
Using this, we get the following version of Theorem \ref{BaPeTh}.

\bigskip

\begin{prop}
\label{asyvart}
For a linear process (\ref{linprodef}) and almost every $\theta\in (0, 2\pi)$, (\ref{norparsumc})  is asymptotically normally distributed under $\mathbb{P}_{\omega}$, for $\mathbb{P}$-almost every $\omega$,  with independent real and imaginary parts, each with mean zero and variance
$$\sigma_{\theta}^2=\frac{|\,||\xi_0||_{_2} f(\theta)|^2}{2},$$
where $f$ is given by (\ref{detpart}).
\end{prop}

By \cite{CuMePe}, p.4075 (Section 4.1) applied to the sequence $(\delta_{1j})_{j\in \mathbb{Z}}$ ($\delta_{ij}$ denotes the Kronecker $\delta-$function) and the fact that $T$ is weakly mixing, the following Law of the Iterated Logarithm holds: for every $t\in (0,2\pi)\setminus\{\pi\}$
\begin{equation}
\label{lilperlin}
\limsup_{n\to \infty}\frac{|\zeta_{-n}(\theta)|}{\sqrt{n \log \log n}}=||\xi_{0}||_{_2}.
\end{equation}

almost surely (note that, for the linear process $(\xi_{n})_{n\in \mathbb{Z}}$, corresponding to convolution with $(\delta_{1j})_{j\in\mathbb{Z}}$, the spectral density with respect to Lebesgue measure is the constant function $||\xi_{0}||_{_2}^{2}/2\pi$).

If $\theta=0$ or $\theta=\pi$, the L.I.L. as stated above holds with $||\xi_{0}||_{_2}$ replaced by $\sqrt{2}\,||\xi_{0}||_{_2}$ (in this case the process $(\zeta_{n}(0))_{n}$ is real-valued). 

The equality (\ref{lilperlin}) clearly implies that $\limsup_{n}\frac{|\zeta_{-n}(\theta)|}{\sqrt{n}}\to \infty$ {\it a.s}. The following lemma states that the divergence occurs ``at comparable speeds'' for every $\theta$. 

\bigskip

\begin{lemma}
\label{lemxij}
Consider the linear process (\ref{linprodef}) and define $\zeta_{-k}$ as in (\ref{defzeta}). Then for every $\lambda\in \mathbb{R}$ and every $0<\eta\leq 1$ there exists an $N\in\mathbb{N}$ satisfying

$$\mathbb{P}\left(\max_{1\leq n\leq N}\frac{|\zeta_{-n}(\theta)|}{\sqrt{n}}> \lambda\right)\geq 1-\eta $$

for all $\theta\in [0,2\pi)$. In particular

$$\mathbb{P}\left(\max_{1\leq n\leq m}\frac{|\zeta_{-n}(\theta)|}{\sqrt{n}}\geq \lambda\right)\geq 1-\eta $$

for all $m\geq N$.
\end{lemma}

{\bf Proof:} Fix $\lambda\in \mathbb{R}$ and $0<\eta\leq 1$. Let\footnote{We work over the invetrval $[0,2\pi]$ (instead of $[0,2\pi)$). This has no effect for the validity of the conclusion and is assumed in order to take advantage of compactness, as will be clear along the proof.} $\theta\in [0,2\pi]$ and $\epsilon>0$ be given and define

$$E_{\epsilon,m}(\theta):=\left[\inf_{|\delta|<\epsilon}\left\{\max_{1\leq n\leq m}\frac{|\zeta_{-n}(\theta+\delta)|}{\sqrt{n}}\right\}>\lambda\right]$$ 
and

$$E_{m}(\theta):=\left[\max_{1\leq n\leq m}\frac{|\zeta_{-n}(\theta)|}{\sqrt{n}}>\lambda\right] .$$

Note that, for fixed  $m$, the sequence of sets $E_{\epsilon,m}(\theta)$ is decreasing with respect to $\epsilon$ ($\epsilon_{1}<\epsilon_{2}$ implies that $E_{\epsilon_{2},m}(\theta)\subset E_{\epsilon_{1},m}(\theta)$), and that the (random) function $\theta\mapsto \max_{1\leq n\leq m}{|\zeta_{-n}(\theta)|}/{\sqrt{n}}$ is continuous for all $m$. In particular

\begin{equation}
\label{decsets}
\bigcup_{\epsilon>0}E_{\epsilon,m}(\theta)=E_{m}(\theta),
\end{equation}

where the union is increasing over $\epsilon$ decreasing to $0$.

By (\ref{lilperlin}), there exists a minimal $N(\theta)$ such that 
$ \mathbb{P}(E_{N({\theta})}(\theta))>1-\eta$ (to see this note that the family $\{E_{k}(\theta)\}_{k\geq 0}$ is increasing with $k$, and its union contains the set $[\limsup_{n}{|\zeta_{-n}(\theta)|}/{\sqrt{n}}>\lambda]$, which has measure $1$ by (\ref{lilperlin})) and therefore, by (\ref{decsets}), there exists an $\epsilon_{\theta}$ such that
\begin{equation}
\label{eqept}
\mathbb{P}(E_{\epsilon_{\theta},N(\theta)}(\theta))>1-\eta \,.
\end{equation}
Now, the family of sets $\{(\theta-\epsilon_{\theta}, \theta+\epsilon_{\theta})\}_{\theta\in [0,2\pi]}$ is an open cover of $[0,2\pi]$, and therefore it admits an open subcover $\{(\theta_{j}-\epsilon_{j}, \theta_{j}+\epsilon_{j})\}_{j=1}^{r}$ (where $\epsilon_{j}:=\epsilon_{\theta_{j}}$).

Let $N=\max\{N(\theta_{1}),\dots,N(\theta_{r})\}$. We claim that, for every $\theta\in [0,2\pi]$
$$\mathbb{P}(E_{N}(\theta))> 1-\eta\,.$$
Indeed, given $\theta\in[0,2\pi]$, with $\theta_{j}-\epsilon_{j}<\theta<\theta_{j}+\epsilon_{j}$,

$$E_{N}(\theta)\supset E_{N(\theta_{j})}(\theta)= \left[\max_{1\leq n \leq N(\theta_{j})}\frac{|\zeta_{-n}(\theta_{j}+\theta-\theta_{j})|}{\sqrt{n}}>\lambda\right]\supset E_{\epsilon_{j},N(\theta_{j})}(\theta_{j})\,,$$

and the conclusion follows from (\ref{eqept}) and the definition of $E_{N}(\theta)$. \qed

\bigskip

\section{The Example}

We now proceed to prove, by an explicit construction, the following proposition:

\bigskip

\begin{prop}
\label{couexafoutra}
There exists a square summable sequence $(a_{n})_{n}$ such that, if $(X_{n})_{n}$ is defined by (\ref{linprodef}) and $S_{n}(\theta)$ is defined by (\ref{nfoutradef}), then
$$\mathbb{P}\left(\limsup_{n\to \infty} \frac{|E_{0}S_{n}(\theta)|}{\sqrt n}=\infty \right)=1$$
for all $\theta\in [0,2\pi)$.

\end{prop}
\bigskip

Before giving the proof we depart from the following observation: if $(n_{k})_{k\geq 0}$ is a strictly increasing sequence of natural numbers and if $(a_{j})_{j}$ is square summable and satisfies $a_{j}=0$ if $j\notin \{n_{k}\}_{k}$  then, using (\ref{exprojsn2})
\begin{equation}
\label{parsnt}
E_{0}S_{n}(\theta)=\sum_{j=0}^{k}e^{in_{j}\theta}a_{n_{j}}(\zeta_{-n_{j}}-\zeta_{-n_{j}+n})(\theta)+ \sum_{j=k+1}^{\infty}e^{in_{j}\theta}a_{n_{j}}(\zeta_{-n_{j}}-\zeta_{-n_{j}+n})(\theta)=:$$
$$A_{k}(n,\theta)+B_{k}(n,\theta)
\end{equation}
so that 

$$\mathbb{P}\left(\max_{n_{k-1}<n\leq n_{k}}\frac{|E_{0} S_{n}(\theta)|}{\sqrt{n}}\geq 2^{k}\right)\geq$$
$$ \mathbb{P} \left(\max_{n_{k-1}<n\leq {n_{k}}}\frac{|A_{k}(n,\theta)|}{\sqrt{n}} \geq 2^{k+1}\right)- \mathbb{P}\left(\max_{n_{k-1}<n\leq {n_{k}}}\frac{|B_{k}(n,\theta)|}{\sqrt{n}}\geq 2^{k}\right)
\geq $$

\begin{equation}
\label{finestsn}
\mathbb{P} \left(\max_{n_{k-1}<n\leq {n_{k}}}\frac{|A_{k}(n,\theta)|}{\sqrt{n}} \geq 2^{k+1}\right)- \mathbb{P}\left(\max_{n_{k-1}<n\leq {n_{k}}}|B_{k}(n,\theta)|\geq 2^{k}\right).
\end{equation}

Now, if $n_{k-1}< n \leq n_{k}$ then, actually

$$A_{k}(n,\theta)= \sum_{j=0}^{k-1}e^{in_{j}\theta}a_{n_{j}}\zeta_{-n_{j}}(\theta)+e^{in_{k}}a_{n_{k}}(\zeta_{-n_{k}}-\zeta_{-n_{k}+n})(\theta).$$
The first summand at the right hand side in this expression is bounded by 
$$\sum_{j=1}^{k-1}\sum_{r=0}^{n_j}|a_{n_{j}}||\xi_{-r}|$$

and therefore there exists $\lambda_{k}>0$ such that

\begin{equation}
\label{firestakn}
\mathbb{P}\left(\left|\sum_{j=0}^{k-1}e^{in_{j}\theta}a_{n_{j}}\zeta_{-n_{j}}(\theta)\right|>\lambda_{k}\right)\leq \left(\frac{1}{2}\right)^{k+2}
\end{equation}

for all $\theta\in [0, 2\pi]$.

All together (\ref{parsnt}), (\ref{finestsn}) and (\ref{firestakn}) give the following result.

\bigskip

\begin{lemma}
\label{lemest}
Let $(n_{k})_{k}$ be a strictly increasing sequence of natural numbers and  let $(a_{j})_{j}$ be a square summable sequence with $a_{j}=0$ for $j\notin \{n_{k}\}_{k}$. Then for every sequence of real numbers $(\lambda_{k})_{k}$ satisfying (\ref{firestakn}) the following inequality holds

$$\mathbb{P}\left(\max_{n_{k-1}< n \leq n_{k}}\frac{|E_{0}S_{n}(\theta)|}{\sqrt{n}}\geq 2^{k}\right)\geq \mathbb{P}\left(\max_{n_{k-1}< n \leq n_{k}}\frac{|(\zeta_{-n_{k}}-\zeta_{-n_{k}+n})(\theta)|}{\sqrt{n}}\geq \frac{\lambda_{k}+2^{k+1}}{a_{n_{k}}} \right)$$
 \begin{equation}
\label{genestpe0}
- \mathbb{P}\left(\max_{n_{k-1}<n\leq {n_{k}}}|B_{k}(n,\theta)|\geq 2^{k}\right)-\left(\frac{1}{2}\right)^{k+2}
\end{equation}
for all $\theta\in [0,2\pi]$.
\end{lemma}

This, together with the previous lemmas, gives the pieces to construct the example stated in Proposition \ref{couexafoutra}.

\bigskip

{\bf Proof of Proposition \ref{couexafoutra}:} 
Following \cite{VoWo}, assume that $||\xi_{0}||_{_2}=1$ and let $(n_{j})_{j\geq 0}$, $(a_{j})_{j\geq 0}$, and $(\lambda_{j})_{j\geq 0}$  be defined inductively as follows:  $n_{0}=1$, $\lambda_{0}=0$, $a_{0}=0$, $a_{1}=\frac{1}{2}$, and given $n_{0},\cdots, n_{k-1}$, $a_{0}, \dots, a_{n_{k-1}}$ and $\lambda_{0},\dots,\lambda_{k-1}$ , let $\lambda_{k}$ be such that
$$\mathbb{P}\left(\left|\sum_{j=1}^{k-1}a_{n_{j}}e^{in_{j}\theta}\zeta_{-n_{j}}(\theta)\right|> \lambda_{k}\right)\leq \left(\frac{1}{2}\right)^{k+2}, $$
 
(see(\ref{firestakn})) and let $n_{k}>n_{k-1}$ be such that

\begin{equation}
\label{almthe}
\mathbb{P}\left(\max_{n_{k-1}< n\leq n_{k}}\frac{|(\zeta_{-n_{k-1}}-\zeta_{-n_{k-1}+n})(\theta)|}{\sqrt{n}}\geq \frac{\lambda_{k}+2^{k+1}}{a_{n_{k-1}}}\right)\geq 1-\left(\frac{1}{2}\right)^{k+1}
\end{equation}

for all $\theta\in[0,2\pi]$. The choice of $n_{k}$ is possible according to Lemma \ref{lemxij} ($|(\zeta_{-n_{k-1}}-\zeta_{-n_{k-1}+n})(\theta)|$ and $|\zeta_{-n}(\theta)|$ have the same distribution). Then define $a_{n_{k}}=\frac{1}{2^{k}\sqrt{n_{k-1}}}$ and $a_{j}=0$ for $n_{k-1}<j<n_{k}-1$.

The sequences $(a_{j})_{j\geq 0}$ and $(\lambda_{k})_{k}$, thus defined, satisfy the hypotheses of Lemma \ref{lemest} and therefore, by the estimates (\ref{genestpe0}) and (\ref{almthe}), 
$$\mathbb{P}\left(\max_{n_{k-1}< n \leq n_{k}}\frac{|E_{0}S_{n}(\theta)|}{\sqrt{n}}\geq 2^{k}\right)\geq 
1-\left(\frac{1}{2}\right)^{k+2}-\mathbb{P}\left(\max_{n_{k-1}<n\leq {n_{k}}}|B_{k}(n,\theta)|\geq 2^{k}\right)$$
for all $\theta\in [0,2\pi]$.

Now we show that, under the present conditions,
\begin{equation}
\label{estbk}
\mathbb{P}\left(\max_{n_{k-1}<n\leq {n_{k}}}|B_{k}(n,\theta)|\geq 2^{k}\right)\leq \left(\frac{1}{2}\right)^{k+2}
\end{equation}
for $k\geq 3$.

Fix $k\geq 3$. First we recall the following (Doob's) maximal inequality for martingales (see \cite{Doob}): if $(M_{n})_{n}$ is a $L^p$ submartingale (namely $||M_{n}||_{p}:=(E|M_{n}|^{p})^{1/p}<\infty$ for all $n$) for some $p>1$ then
\begin{equation}
\label{doomarine}
||\sup_{k\leq n}M_{k}||_{p}\leq \frac{p}{p-1}||M_{n}||_{p} \,.
\end{equation}
Now, for fixed $\theta$, $(|\zeta_{-n}(\theta)|)_{n\geq 0}$ is an $L^2$ submartingale (with respect to $(\mathcal{G}_{n})_{n}$, where $\mathcal{G}_{k}=\sigma((\xi_{-j})_{j\leq k-1})$) and therefore, by Doob's maximal inequality (\ref{doomarine}):

$$E\left(\max_{k\leq n}\left|\zeta_{-k}(\theta)\right|\right)\leq \left|\left|\max_{k\leq n}\left|\zeta_{-k}(\theta)\right|\right|\right|_{2}\leq 2\,||\zeta_{-n}(\theta)||_{_2}\leq 2\,\sqrt{n}.$$

This gives

$$E\left(\max_{n_{k-1}<n\leq n_{k}}\left|B_{k}(n,\theta)\right|\right)\leq \sum_{j=k+1}^\infty a_{n_{j}}E\left(\max_{k\leq n_{k}-n_{k-1}}\left|\zeta_{-k}(\theta)\right|\right)\leq$$

$$\sum_{j=k+1}^\infty \frac{1}{2^{j-1}}\sqrt{\frac{n_{k}-n_{k-1}}{n_{j-1}}}\leq \frac{1}{2^{k-1}}, $$
  and therefore, by Markov's inequality
  
  $$\mathbb{P}\left(\max_{n_{k-1}<n\leq n_{k}}\left|B_{k}(n,\theta)\right|\geq 2^{k}\right)\leq \frac{1}{2^{^{2k-1}}}\leq \left(\frac{1}{2}\right)^{k+2}$$
  
 as claimed.
 
 The proof is finished as follows: a combination of (\ref{genestpe0}), (\ref{almthe}) and (\ref{estbk}) gives, under the present choices of $(a_{k})_{k}$ and $(n_{k})_{k}$, that
 
 $$\mathbb{P}\left(\max_{n_{k-1}< n \leq n_{k}}\frac{|E_{0}S_{n}(\theta)|}{\sqrt{n}}< 2^{k}\right)\leq \left(\frac{1}{2}\right)^{k+1}$$
 so that, by the first Borel-Cantelli Lemma
  $$\max_{n_{k-1}< n \leq n_{k}}\frac{|E_{0}S_{n}(\theta)|}{\sqrt{n}}\geq 2^{k} \mbox{   \it except for finitely many $k$'s,}$$
$\mathbb{P}-$a.s. This clearly implies that $\limsup_{n}{|E_{0}S_{n}(\theta)|}/{\sqrt{n}}=\infty$ {$\mathbb{P}-$a.s.}  \qed
  
\section*{Appendix: convergence of Types}

A distribution function $F$ is {\it non-degenerate} if it is not a Heaviside function (this is, if it is not the indicator function of some interval $[a,+\infty)$). We recall the following {\it Convergence of Types} theorem (\cite{Bil}, Th 14.2).

\medskip

\begin{lemma}
\label{BilConTyp}
Let $F_{n}$, $F$ and $G$ be distribution functions,  and $a_{n},u_{n},b_{n}, v_{n}$ be constants with $a_{n}>0,\,u_{n}>0$. If $F$, $G$ are nondegenerate, $F_{n}(a_{n}x+b_{n})\Rightarrow F(x)$, and $F_{n}(u_{n}x+v_{n})\Rightarrow G(x)$ then there exist $a=\lim_{n}a_{n}/u_{n}$, $b=\lim_{n}(b_{n}-v_{n})/u_{n}$ and $G(x)=F(ax+b)$.  
\end{lemma}

Note that, necessarily, $a>0$ (as otherwise $G$ would be constant).

We will translate this statement to a statement about convergence of stochastic processes (with a restricted choice of $u_{n},v_{n}$, see Proposition \ref{quecon} below). To begin, we remind the following elementary facts, here $\to_{P}$ denotes convergence in probanbility.

\begin{enumerate}

\item If $a$ is constant then $U_{n}\Rightarrow a$ if and only if $U_{n}\to_{P} a$. 

\item If $U_{n}\Rightarrow W$ and $V_{n}\to_{P} 0$  then $U_{n}+V_{n}\Rightarrow W$.

\item If $(a_{n})_{n}$ is a sequence of constant functions then $a_{n}\Rightarrow A$ if and only if $a=\lim_{n}a_{n}$ exists (and therefore $A=a$ a.s.).

\end{enumerate}

\bigskip

\begin{lemma}
\label{lemcondeg}
If  $Y_{k}\Rightarrow Y$ and $\{c_{k}\}_{k}\subset \mathbb{R}$ are such that $Y_{k}+c_{k}\Rightarrow 0$, then $Y=-\lim_{k}c_{k}$. In particular, $Y$ is a constant function. 
\end{lemma}

{\bf Proof:} Note that $c_{k}=-Y_{k}+(Y_{k}+c_{k})\Rightarrow -Y$ because $Y_{k}+c_{k}\Rightarrow 0$. Now use 2. and 3. above. \qed

\bigskip

\begin{cor}
\label{degcas}
If $X$ is not constant, $X_{n}\Rightarrow X$, and $a_{n}$, $b_{n}$ are such that $a_{n}X_{n}+b_{n}\Rightarrow 0$, then $a_{n}\to 0$ and $b_{n}\to 0$.
\end{cor}

{\bf Proof:} If $0<a:=\limsup_{n}a_{n}\leq \infty$ and $a_{n_{k}}\to_{k\to \infty} a$ with $a_{n_{k}}>0$, then applying Lemma \ref{lemcondeg} with $Y_{k}=X_{n_{k}}$ and $c_{k}=b_{n_{k}}/a_{n_{k}}$ we conclude that $X$ is constant. This proves that, necessarily, $\limsup_{n}a_{n}\leq 0$. A similar argument shows that $\liminf_{n}a_{n}\geq 0$, and therefore $\lim_{n}a_{n}=0$.

The fact that $b_{n}\to 0$ follows from here applying Lemma \ref{lemcondeg} again, because $a_{n}X_{n}\Rightarrow 0$ . \qed  

\medskip

These results give rise to the following proposition

\bigskip

\begin{prop}
\label{quecon}
If $X$ is not constant, $X_{n}\Rightarrow X$ and $a_{n}>0$, $b_{n}$ are such that $a_{n}X_{n}+b_{n}\Rightarrow Y$, then there exists $a=\lim_{n}a_{n}$, $b=\lim_{n}b_{n}$ and, therefore, $Y=aX+b$ (in distribution).
\end{prop}

{\bf Proof:} If $Y$ is constant then, from $a_{n}X_{n}+b_{n}-Y\Rightarrow 0$ (see 1. above) it follows, via Corollary \ref{degcas}, that $\lim_{n}a_{n}=0$ and $\lim_{n}b_{n}=Y$.

If $Y$ is not constant we apply Lemma \ref{BilConTyp} with $F_{n}$, $F$, and $G$ the distribution functions of $X_{n}$, $X$ and $Y$ respectively, and with $u_{n}=1$, $v_{n}=0$.\qed

\medskip
 
{\bf Remark:} Taking $X_{n}=1$ (the constant function), $a_{n}=n$, and $b_{n}=-n$, we see that the given restriction on $X$ (to be non constant) is necessary.

\medskip

We finish this appendix with a lemma covering the asymptotically degenerate case.

\medskip

\begin{lemma}
\label{degcascon}
If $Y$ is constant, $Y_{k}\Rightarrow Y$, and $Y_{k}+c_{k}\Rightarrow Z$ then $c=\lim_{k}c_{k}$ exists and $Z=Y+c$. 
\end{lemma}

{\bf Proof:} Use $Y+c_{k}=(Y-Y_{k})+Y_{k}+c_{k}\Rightarrow Z$ by 2. and 3. above.\qed

\bigskip
  
{\bf Acknowledgements}

The author wants to thank Magda Peligrad for proposing the question of this note and for our discussions related to it. Valuable comments and suggestions are also due to Dalibor Voln\'{y}. The author was supported by the NSF Grant DMS-1208237.

\end{document}